\documentclass[11pt, a4paper,twoside]{article}
\usepackage{enumerate}
\usepackage{amsmath}
\usepackage{fix-cm}
\usepackage{pb-diagram}
\usepackage[english]{babel}
\usepackage{amssymb,latexsym}
\usepackage{latexsym}
\def \bui#1#2{\mathrel{\mathop{\kern 0pt#1}\limits^{#2}}}

\let\<\langle
\let\>\rangle
\newcommand{\R}{{\mathbb R}}

\newtheorem{example}{Examples}[section]
\newtheorem{thm}{Theorem}[section]
\newtheorem{lemma}[thm]{Lemma}
\newtheorem{prop}[thm]{Proposition}

\newtheorem{cor}[thm]{Corollary}
\newtheorem{remark}[thm]{Remark}
\newtheorem{remarks}[thm]{Remarks}
\newtheorem{definition}[thm]{Definition}
\newtheorem{notation}[thm]{Notation}
\newtheorem{exabout:ample}[thm]{Example}

\title{The Bochner formula for Riemannian flows}
\author{Fida El Chami\footnote{Lebanese University, Faculty of Sciences II, Department of Mathematics, P.O. Box 90656 Fanar-Matn, Lebanon,
E-mail: \texttt{fchami@ul.edu.lb}},\, Georges Habib\footnote{Lebanese University, Faculty of Sciences II, Department of Mathematics, P.O. Box 90656 Fanar-Matn, Lebanon,
E-mail: \texttt{ghabib@ul.edu.lb}}}

\begin{document}
\date{}
\maketitle
\begin{abstract} \noindent In this paper, we consider a Riemannian manifold $(M,g)$ endowed with a Riemannian flow and we study the curvature term in the Bochner-Weitzenb\"ock formula of the basic Laplacian on $M.$ We prove that this term splits into two parts. The first part depends mainly on the curvature operator of the underlying manifold $M$ and the second part is expressed in terms of the O'Neill tensor of the flow. After getting a lower bound for this term, depending on these two parts, we establish an eigenvalue estimate of the basic Laplacian on basic forms. We then discuss the limiting case of the estimate and prove that when equality occurs, the manifold $M$ is a local product. This paper follows mainly the approach in \cite{S}. 
\end{abstract}

\section{Introduction} Given a Riemannian manifold $(M^n,g)$ and a form $\omega$ on $M$ of degree $p,$ the Laplacian $\Delta$ of $\omega$ is related to the curvature operator on $M$ through the Bochner-Weitzenb\"ock formula, namely 
$$\Delta\omega=\nabla^*\nabla \omega+\mathcal{B}^{[p]}\omega,$$ 
where $\mathcal{B}^{[p]},$ usually called the {\it Bochner operator}, is the symmetric endomorphism of the bundle of $p$-forms $\Lambda^p(M)$ given by $\mathcal{B}^{[p]}=\mathop\sum_{i,j=1}^n e_j\wedge (e_i\lrcorner R^M(e_j,e_i)).$ Here $R^M$ is the curvature operator on $M$ defined by convention $R^M(X,Y)=\nabla^M_{[X,Y]}-[\nabla^M_X,\nabla^M_Y]$ and $\{e_i\}_{i=1,\cdots,n}$ denotes a local orthonormal frame of $TM.$ In all the paper, we identify vector fields with their corresponding $1$-forms through the usual musical isomorphisms. 

\noindent It is clear that the Bochner-Weitzenb\"ock formula is a useful tool to estimate the eigenvalues of the Laplacian ($M$ is assumed to be compact and connected in this case), since any lower bound of the Bochner operator provides a lower bound of the eigenvalues. For example, when $p=1,$ A. Lichnerowicz \cite{L} proved that if $\mathcal{B}^{[1]}$ (which corresponds to the Ricci tensor of the manifold) is greater than some positive number $k,$ the first positive eigenvalue is greater than $k\frac{n}{n-1}.$ This inequality was later characterized by M. Obata in \cite{Ob} who states that equality occurs if and only if the manifold is isometric to a round sphere.  

\noindent Another estimate of the Bochner operator were obtained by Gallot and Meyer in \cite{GM} when $p=1,\cdots,n-1.$ Indeed, they showed that if the curvature operator of $M$ has a lower bound $k,$ then $\mathcal{B}^{[p]}$ is always greater than $p(n-p)k.$ This inequality has led to the following rigidity result \cite[Prop. 2.9]{GM}: when the lower bound $k$ is strictly positive, then all the cohomology groups $H^p(M)$ vanish which mainly means that the manifold $M$ has the same cohomology as the round sphere. Moreover, based on the same inequality, they proved the following estimates for the first eigenvalue of the Laplacian restricted to closed forms $\lambda'_{1,p}$ and to co-closed forms $\lambda''_{1,p}$, namely 
\begin{equation}\label{eq:gallotmeyer}
\lambda'_{1,p}\geq  kp(n-p+1) \quad \text{and}\quad \lambda''_{1,p}\geq  k(p+1)(n-p).
\end{equation}  
Here $k$ is assumed to be strictly positive. Besides the round sphere of curvature $k,$ the authors provided examples of hypersurfaces in the complex projective space where the equality in \eqref{eq:gallotmeyer} is attained \cite[Prop. 8.1]{GM}. 

\noindent In \cite{S}, the author used a new technique to bound the Bochner operator on submanifolds. In fact, on a given Riemannian manifold $M$ of dimension $n$ and a submanifold $\Sigma,$ he expressed the curvature operator on $\Sigma$ in terms of the one on $M$ and the second fundamental form  of the immersion through the Gauss formula. Namely, he showed that the term $\mathcal{B}^{[p]},$ acting on $p$-forms of $\Sigma,$ can be splitted into two parts: the restriction part $\mathcal{B}_{{\rm res}}^{[p]}$ that mainly depends on the ambient manifold $M$ and the exterior part $\mathcal{B}_{{\rm ext}}^{[p]}$ that is determined by the Weingarten tensor $S$ \cite[Thm. 1]{S}. The proof is based on the expression of the Bochner operator $\mathcal{B}^{[p]}$ in terms of the curvature of the underlying manifold $\Sigma$ through the Clifford multiplication used in \cite{P}. More precisely, for hypersurfaces, he proved that the following inequality 
$$\mathcal{B}^{[p]}\geq p(n-p)(\gamma_M+\beta_p(\Sigma)),$$
holds, where $\gamma_M$ is a lower bound of the curvature operator of $M$ and $\beta_p(\Sigma)$ is the lowest eigenvalue of the operator $T^{[p]}=({\rm tr}\, S)S^{[p]}-S^{[p]}\circ S^{[p]}.$ The operator $S^{[p]}$ is some canonical extension of $S$ to $p$-forms on $\Sigma.$ After estimating the eigenvalues of the tensor $T^{[p]}$ in terms of different geometric quantities on $\Sigma$ such as the mean curvature and the norm of the Weingarten tensor, he deduced several rigidity results, among them  the de Rham cohomology groups of $\Sigma,$ certain Clifford torus and immersions of K\"ahler manifolds... In the same spirit and with the use of the Bochner-Weitzenb\"ock formula, he found a sharp estimate for the eigenvalues of the Laplacian on $\Sigma$ that involves geometric data of the immersion. We note that this eigenvalue estimate has been later generalized to all codimensions in \cite{CS}.

\noindent In this paper, we study the Bochner operator for Riemannian flows (see Section \ref{sec:1} for the definition). These are the global geometric aspects of Riemannian submersions. Roughly speaking, a Riemannian flow on a given Riemannian manifold $(M,g)$ is determined by a unit vector field $\xi$ on $M$ such that the Lie derivative of the metric $g$ vanishes along $\xi$ when one restricts to vector fields orthogonal to $\xi.$ Examples of Riemannian flows are provided by Killing vector fields, Sasakian manifolds...We notice here that the integral curves of $\xi,$ called {\it the leaves}, are the fibers of local Riemannian submersions that map to a Riemannian manifold which detects the {\it transverse geometry} of the flow. When looking to the structure of the normal bundle $Q=\xi^\perp$ of the flow, we require objects to be {\it basic} which means that they just depend on the transverse variables. In this spirit, a transverse Bochner-Weitzenb\"ock formula carries over for the {\it basic Laplacian} (see Equation \eqref{bochner}); this allows to study the geometric and analytic properties of the flow, such as the basic cohomology groups. 

\noindent Following the approach of A. Savo in \cite{S}, we consider in this work the Bochner operator in the transverse Bochner-Weitzenb\"ock formula and we aim to express it in terms of the geometric data of the flow. We prove with the help of the O'Neill formulas \cite{O} that, as for submanifolds, the Bochner operator splits into a restriction part and an exterior part (see Equation \eqref{eq:decomposition}) where the first part depends on the geometry of the ambient manifold while the second part involves the O'Neill tensor. Using this expression, we deduce a lower bound of this operator in Corollary \ref{curvatureestimate3} which allows to get vanishing results on the basic cohomology groups (see Corollary \ref{corrigidity}). Also in Theorem \ref{estimateeigenvalues}, we establish a sharp estimate for the first eigenvalue $\lambda_{1,p}$ of the basic Laplacian restricted to $p$-forms ($1\leq p\leq [\frac{q}{2}]$ with $q$ is the codimension of the flow). Namely, we show 
$$\lambda_{1,p} \ge p(q-p+1)(\gamma_M +\beta_M^1),$$
where $\gamma_M$ is a lower bound of the curvature operator on $M$ restricted to $Q$ and $\beta_M^1$ is the lowest eigenvalue of the symmetric tensor $h^2$ ($h$ denotes the O'Neill tensor). When equality occurs in the above estimate, we show that the O'Neill tensor vanishes and the manifold $M$ is then isometric to the quotient of $\mathbb{R}\times \Sigma$ by some subgroup $\Gamma,$ where $\Sigma$ is a compact simply connected manifold of positive curvature. 

\noindent The paper is organized as follows. In Section \ref{sec:1}, we review the definitions of foliations and the basic Laplacian. We also state an eigenvalue estimate for the basic Laplacian that involves a lower bound of the Bochner operator (see Proposition \ref{estimate}). In Section  \ref{sec:clifford}, we adapt the way of writing the Bochner operator in terms of Clifford multiplication used in \cite{P} to the set-up of foliations. We then prove a rigidity result of the basic cohomology groups stating that they all vanish when the transverse curvature operator is positive (see Proposition \ref{vanishing}). The main results are stated and proven in Section \ref{sec:flows} where the case of Riemannian flow is considered. The last section is devoted to a well-known general results on foliations that we use in our study.

{\bf Acknowledgment:} We would like to thank Ola Makhoul, Nicolas Ginoux and Ken Richardson for many helpful discussions during the preparation of this paper. This project is supported by a grant from the Lebanese University.

\section{Preliminaries}\label{sec:1} 

\noindent In this section, we recall the main definitions on Riemannian foliations and some known results that can be found in \cite{T}.

\noindent Let $(M^n,g, \mathcal{F})$ be a Riemannian manifold of
dimension $n$ endowed with a Riemannian foliation $\mathcal{F}$ of codimension $q$. We assume, throughout this paper, that the metric $g$ is bundle-like \cite{T}. That means, $\mathcal{F}$ is given by an integrable subbundle $L$ of $TM$ and the metric $g$ satisfies the condition $\mathcal{L}_Xg|_Q=0$ on the normal vector bundle $Q=TM/L,$ for all $X\in \Gamma(L).$ Here $\mathcal{L}$ denotes the Lie derivative. In this case, the tangent bundle of $M$ decomposes orthogonally into $L$ and $Q.$ We equip the normal bundle $Q$ with the transverse Levi-Civita connection $\nabla$ \cite{T}. It is a standard fact that the curvature associated to $\nabla$ vanishes along the leaves and therefore curvature data on $Q$ are defined along orthogonal directions. Recall also that a basic form $\omega$ is a differential form on $M$ that does uniquely depend on the transverse variables, in other words, $\omega$ satisfies $X\lrcorner\omega=0$ and $X\lrcorner d\omega=0,$ for all $X\in \Gamma(L)$. These basic forms are preserved by the exterior derivative and are used to define the basic Laplacian $\Delta_b=d_b\delta_b+\delta_bd_b$ where $d_b$ is the restriction of the exterior differential $d$ to basic forms and $\delta_b$ is its $L^2$-adjoint. The basic Laplacian yields the basic Hodge theory that can be used to compute the basic cohomology groups $H^p_b(\mathcal{F})=\frac{{\rm ker}\, d_{p}}{{\rm image}\,d_{p-1}}.$ In the study of the basic Poincar\'e duality (which fails to hold for the basic Laplacian), the authors in \cite{HR} introduce a new cohomology ${\widetilde H}_b(\mathcal{F})$ that uses the {\it twisted exterior derivative} $\widetilde{d}_b:=d_b-\frac{1}{2}\kappa_b\wedge,$ where $\kappa_b$ is the basic component of the mean curvature field $\kappa$ of the foliation \cite{RP}. They prove that the associated twisted Laplacian $\widetilde{\Delta}_b:=\widetilde{d}_b\widetilde{\delta}_b+\widetilde{\delta}_b\widetilde{d}_b$ commutes with the basic Hodge operator and therefore the Poincar\'e duality carries on for those twisted cohomology groups. Also, they state a Bochner-Weitzenb\"ock formula for $\widetilde{\Delta}_b$ which allows to generalize several rigidity results on the usual basic cohomology. Namely, on basic $p$-forms, the formula is \cite[Prop. 6.7]{HR} 
\begin{equation}\label{bochner}
\widetilde{\Delta}_b=\nabla^*\nabla+\mathcal{B}^{[p]}+\frac{1}{4}|\kappa_b|^2,
\end{equation}
where $\nabla^*\nabla:=-\sum_{i=1}^q\nabla_{e_i,e_i}+\nabla_{\kappa_b}$ and $\mathcal{B}^{[p]}=\mathop\sum_{i,j=1}^q e_j\wedge (e_i\lrcorner R(e_j,e_i))$ with $R(X,Y)=\nabla_{[X,Y]}-[\nabla_X,\nabla_Y]$ is the transversal curvature operator, $\{e_i\}_{i=1,\cdots,q}$ is a local orthonormal frame of $Q$. Here the basic component of the mean curvature $\kappa_b$ is assumed to be a harmonic $1$-form. As the spectrum of $\widetilde{\Delta}_b$ (as well as the dimensions of ${\widetilde H}_b(\mathcal{F})$) remains invariant for any choice of the bundle-like metric \cite{HR}, one can state the following, as in \cite[Prop. 3]{S}, 

\begin{prop} \label{estimate} Let $(M,g,\mathcal{F})$ be a compact Riemannian manifold endowed with a Riemannian foliation $\mathcal{F}$ of codimension $q$ and a bundle-like metric $g.$ Let $p$ be an integer number such that $1\leq p\leq q.$ 
\begin{enumerate} 
\item [1)] If $\mathcal{B}^{[p]}\geq 0$ and $\kappa_b$ is a basic-harmonic one form, then any basic harmonic $p$-form is transversally parallel. If the strict inequality $\mathcal{B}^{[p]}>0$ holds, then $H_b^p(\mathcal{F})=0.$ 
\item [2)] If the foliation is minimal and $\mathcal{B}^{[p]}\geq p(q-p)\Lambda$ for some $\Lambda>0,$ then the first eigenvalue $\lambda_{1,p}$ of the basic Laplacian satisfies
\begin{equation}\label{eq:estimatebasic}
\lambda_{1,p}\geq p(q-p+1)\Lambda,
\end{equation}
where $p$ is chosen such that $1\leq p\leq \frac{q}{2}.$ 
\end{enumerate}
\end{prop} 

{\noindent \it Proof.} The proof of the point $1)$ is a direct consequence of the Bochner-Weitzenb\"ock formula. Indeed, take any basic harmonic $p$-form $\omega,$ that is $d_b\omega=\delta_b \omega=0,$ one can easily see that $|{\widetilde d}_b\omega|^2+|{\widetilde \delta}_b\omega|^2=\frac{1}{4}|\kappa_b|^2|\omega|^2.$ Hence, applying Equation \eqref{bochner} to $\omega$ and taking the scalar product the same form, one gets after integrating over $M$
$$\frac{1}{4}\int_M |\kappa|^2|\omega|^2 dv_g=\int_M|\nabla\omega|^2 dv_g+\int_M\langle \mathcal{B}^{[p]}\omega,\omega\rangle dv_g+\frac{1}{4}\int_M |\kappa|^2|\omega|^2 dv_g\geq \frac{1}{4}\int_M |\kappa|^2|\omega|^2 dv_g,$$
which allows to deduce the first statement. Now, if $\mathcal{B}^{[p]}>0$ then it is clear that any basic harmonic $p$-form vanishes. By \cite{D} and \cite[Thm 6.2]{M}, one can always change the bundle-like metric into another bundle-like metric (with the same transverse metric) so that the basic component of the mean curvature $\kappa_b$ is a basic harmonic $1$-form with respect to the new metric. Therefore, we can work with such a metric keeping the same condition on $\mathcal{B}^{[p]}.$ Hence the assumption on the mean curvature can be dropped off and we deduce the statement 2) as the basic cohomology is independent of the choice of the bundle-like metric. The proof of the point $2)$ follows the same way as in \cite{GM} by proving that $\int_M |\nabla \omega|^2 dv_g\geq \frac{\lambda_{1,p}}{q-p+1}\int_M |\omega|^2 dv_g,$ which finishes the proof. 
\hfill$\square$ 

{\noindent\bf Remark.} We point out that when the equality case in \eqref{eq:estimatebasic} is attained, the associated eigenform is a {\it basic conformal Killing form} \cite{Se, JR} which is either closed or of degree $p=\frac{q}{2}$ (that is, $q$ should be even). Recall here that a basic conformal Killing form $\omega$ is a basic form that satisfies, for all $X\in \Gamma(Q),$ the equation 
$$\nabla_X\omega=\frac{1}{p+1}X\lrcorner d\omega-\frac{1}{q-p+1}X\wedge \delta_T\omega,$$ 
where $\delta_T=\delta_b-\kappa_b\lrcorner.$ 

\section{Clifford multiplication on basic forms} \label{sec:clifford}

\noindent In this section, we will review the approach of \cite[Sect. 4]{P} to write the curvature term in the Bochner-Weitzenb\"ock formula in terms of the Clifford multiplication. We also refer to \cite{S} for more details. 

\noindent Let $(M,g,\mathcal{F})$ be a Riemannian manifold endowed with a Riemannian foliation $\mathcal{F}$ and let $Q$ be the normal bundle of codimension $q$. For $X\in \Gamma(Q)$ and $\omega$ a $p$-form on $Q$, the Clifford multiplication of $X$ with $\omega$ is defined as 
\begin{equation}\label{cliffordmul}
X\cdot\omega= X\wedge\omega-X\lrcorner \omega \quad\text{and}\quad \omega\cdot X= (-1)^p (X\wedge\omega+X\lrcorner \omega).
\end{equation}
A direct consequence of the definitions says that for any two sections $X$ and $Y$ on $Q$, the following relation 
$$X\cdot Y+Y\cdot X=-2g(X,Y)$$
holds. Given any two forms $\omega_1$ and $\omega_2,$ one can extend the definition \eqref{cliffordmul} to the Clifford multiplication between $\omega_1$ and $\omega_2$ as follows: write locally $\omega_1=\mathop\sum_{i_1\leq \cdots \leq i_p} \alpha_{i_1\cdots i_p} e_{i_1}\wedge \cdots \wedge e_{i_p}$ in any orthonormal frame $\{e_1,\cdots,e_q\}$ of $Q$ and define 
$$\omega_1\cdot \omega_2=\sum_{i_1\leq \cdots \leq i_p} \alpha_{i_1\cdots i_p} e_{i_1}\cdot \cdots e_{i_p}\cdot \omega_2.$$
The Lie bracket is then defined as $[\omega_1,\omega_2]=\omega_1\cdot\omega_2-\omega_2\cdot\omega_1.$ For a $2$-form $\Psi$ and a $p$-form $\omega$,  the Lie bracket between $\Psi$ and $\omega$ can be expressed explicitly as 

\begin{lemma} \label{bracketcomputation1} Let $\Psi$ be a $2$-form and let $\omega$ be a $p$-form. One has 
$$[\Psi,\omega]=2\sum_{i=1}^q (e_i\lrcorner \psi)\wedge (e_i\lrcorner \omega),$$
where $\{e_1,\cdots,e_q\}$ is an orthonormal frame of $Q.$ In particular the degree of $[\Psi,\omega]$ is the same as the form $\omega.$
\end{lemma}

{\noindent \it Proof.} The proof relies mainly on the use of Equations \eqref{cliffordmul} and the fact that $X\cdot\omega=(-1)^p\omega\cdot X-2X\lrcorner \omega.$ Indeed, if we write $\Psi=\sum_{i<j} \Psi_{ij}e_i\wedge e_j,$ we compute 
\begin{eqnarray*} 
\Psi\cdot\omega&=&\sum_{i<j} \Psi_{ij}e_i\cdot e_j\cdot\omega=\sum_{i<j} \Psi_{ij}e_i\cdot ((-1)^p\omega\cdot e_j-2e_j\lrcorner \omega)\\
&=&\sum_{i<j} \Psi_{ij} \left(\omega\cdot e_i\cdot e_j-2(-1)^p(e_i\lrcorner \omega) \cdot e_j-2e_i\cdot (e_j\lrcorner \omega) \right)\\
&=& \omega\cdot\Psi-2(-1)^p\sum_{i<j} \Psi_{ij}(e_i\lrcorner \omega) \cdot e_j-2(-1)^{p-1}\sum_{i<j} \Psi_{ij}(e_j\lrcorner \omega)\cdot e_i+4\sum_{i<j} \Psi_{ij}e_i\lrcorner (e_j\lrcorner \omega)\\
&=&\omega\cdot\Psi-2(-1)^p\sum_{i,j} \Psi_{ij}(e_i\lrcorner \omega) \cdot e_j+2\sum_{i,j} \Psi_{ij}e_i\lrcorner (e_j\lrcorner \omega)\\
&=& \omega\cdot\Psi+2\sum_{i,j} \Psi_{ij}(e_j\wedge (e_i\lrcorner \omega)+e_j\lrcorner (e_i\lrcorner \omega)) +2\sum_{i,j} \Psi_{ij}e_i\lrcorner (e_j\lrcorner \omega).
\end{eqnarray*}
Finally, we deduce that  $[\Psi,\omega]=2\sum_{i,j} \Psi_{ij}e_j\wedge (e_i\lrcorner \omega)$ which finishes the proof of the lemma. 
\hfill$\square$

\noindent Another useful property of the Lie bracket that will be used later in this paper. 

\begin{lemma} \label{bracketcomputation} Let $\Psi$ be a $2$-form and let $\omega$ be a $p$-form. Then we have
$$[\Psi,X\wedge \omega]=X\cdot[\Psi,\omega]+2(X\lrcorner\Psi)\cdot \omega+[\Psi,X\lrcorner \omega],$$ 
for any $X\in \Gamma(Q).$
\end{lemma}

{\noindent \it Proof.} Using the definition of the Lie bracket, we write 
\begin{eqnarray*} 
[\Psi,X\wedge \omega]&=&\Psi\cdot(X\wedge\omega)-(X\wedge\omega)\cdot \Psi\\
&=&\Psi\cdot(X\cdot\omega+X\lrcorner\omega)-(X\cdot\omega+X\lrcorner\omega)\cdot \Psi\\
&=&X\cdot\Psi\cdot \omega+2(X\lrcorner\Psi)\cdot\omega+\Psi\cdot(X\lrcorner\omega)-X\cdot\omega\cdot \Psi-(X\lrcorner\omega)\cdot \Psi\\
&=& X\cdot[\Psi,\omega]+2(X\lrcorner\Psi)\cdot \omega+[\Psi,X\lrcorner \omega].
\end{eqnarray*}
The proof of the lemma is then finished. 
\hfill$\square$

\noindent Next, we recall the definition of the basic Dirac operator restricted to basic forms \cite{GK}. Given any orthonormal frame $\{e_i\}_{i=1,\cdots,q}$ of $\Gamma(Q),$ the basic Dirac operator is defined as 
$$D_b=\sum_{i=1}^q e_i\cdot\nabla_{e_i}-\frac{1}{2}\kappa_b\cdot,$$ 
where $\kappa_b$ is as usual the projection of the mean curvature. It is easy to see that $D_b=\widetilde{d}_b+\widetilde{\delta}_b$ and that $D_b^2={\widetilde\Delta}_b.$ As in \cite[Thm. 50]{P}, one can show that (see also \cite[Prop. 1.3.5]{Hab}) 
$$D_b^2\omega =\nabla^*\nabla\omega -\frac{1}{2}\sum_{i,j=1}^q e_i\cdot e_j\cdot R(e_i,e_j)\omega +\frac{1}{4}|\kappa_b|^2\omega ,$$
and 
$$D_b^2\omega=\nabla^*\nabla \omega +\frac{1}{2}\sum_{i,j=1}^q  R(e_i,e_j)\omega\cdot e_i\cdot e_j+\frac{1}{4}|\kappa_b|^2\omega.$$
Now by adding these two equations and dividing by $2$, we deduce after comparing with Equation \eqref{bochner} that 
$$\mathcal{B}^{[p]}\omega=\frac{1}{4}[R(e_i,e_j)\omega, e_i\cdot e_j].$$ 
Following the same lines of the proof of \cite[Thm. 17]{S}, one can say that 
\begin{equation}\label{curvature}
\langle \mathcal{B}^{[p]} \omega , \varphi \rangle = \dfrac{1}{4} \sum_{r,s=1}^{\binom{q}{2}} \langle R \psi_r , \psi_s \rangle \langle [\hat{\psi}_r, \omega],  [\hat{\psi}_s, \varphi]    \rangle ,
\end{equation}
where $\{\psi_r\}_{r=1, \dots, \binom{q}{2}}$ is any orthonormal frame of $\wedge^2 Q$ and that  $\{\hat{\psi}_r\}_{r=1, \dots, \binom{q}{2}}$ its dual basis. Here the curvature $R:\Lambda^2Q \to \Lambda^2 Q$ is viewed as a symmetric operator by $\langle R(X\wedge Y),Z\wedge W\rangle=g(R(X,Y)Z,W)$ for all $X,Y,Z,W\in \Gamma(Q).$ 

\noindent As in \cite[Thm. 51]{P}, we deduce the following result (see also \cite[Cor. D]{ORT} for a different proof)

\begin{prop}\label{vanishing}
Let $(M,g,\mathcal{F})$ be a compact Riemannian manifold endowed with a Riemannian foliation of codimension $q$. 
\begin{enumerate}
\item If the transversal curvature operator is nonnegative and $\kappa_b$ is basic-harmonic, then any basic harmonic form is transversally parallel.
\item If the transversal curvature operator is positive, then $H^p_b(\mathcal{F})=0$ for all $p\in \{1,\cdots,q-1\}.$
\end{enumerate}
\end{prop}

\section{Riemannian flows} \label{sec:flows}

\noindent In this section, we will consider a Riemannian flow, that is a Riemannian foliation of $1$-dimensional leaves given by a unit vector field. As mentioned in the introduction, we will prove throughout this section that the curvature operator of the normal bundle splits into two parts. The first part, that we call restriction part, depends mainly on the curvature operator of the underlying manifold and the second part, that we call exterior part, is expressed in terms of the O'Neill tensor of the flow. 

\noindent Let $(M,g,\xi)$ be a Riemannian manifold endowed with a Riemannian flow given by a unit vector field $\xi.$ Recall the condition on the metric that $\mathcal{L}_\xi g|_{\xi^\perp}=0$ which means that the tensor $h=\nabla^M\xi$, called the O'Neill tensor,  is a skew-symmetric endomorphism on $\Gamma(Q).$ From the relation $g(h(X),Y)=-\frac{1}{2}g([X,Y],\xi)$, one can characterize the integrability of the normal bundle of a Riemannian flow by the vanishing of the O'Neill tensor \cite{O}. Moreover, when the O'Neill tensor and the mean curvature $\kappa:=\nabla^M_\xi\xi$ both vanish, the manifold $M$ is isometric to a local product. Also, one can easily check by a straightforward computation that when the mean curvature $\kappa$ is a basic one form, the endomorphism $h$ is a basic tensor, that is, $\nabla_\xi h=0.$ Recall here that $\nabla$ is the transversal Levi-Civita connection extended to tensors. Based on this fact, the curvature $R^M$ restricted to sections of the form $\xi\wedge X$ for $X\in \Gamma(Q)$ can be expressed as follows 

\begin{lemma} \label{computationcurvature} On a Riemannian manifold $(M^n,g,\xi)$ endowed with a Riemannian flow with basic mean curvature $\kappa$, we have  that 
$$R^M(\xi,X)\xi=-h^2(X)+g(\kappa,h(X))\xi+\nabla^M_X\kappa-g(\kappa,X)\kappa,$$ 
for any $X\in \Gamma(Q).$ In particular, for minimal Riemannian flow, the matrix of $R^M$ in the orthonormal frame $\{\xi\wedge e_i\}_{i=1,\cdots,n-1}$ is the same as $-h^2.$ 
\end{lemma}

\noindent{\it Proof.} Let $X$ be any foliated vector field, that is $\nabla_\xi X=0.$ The curvature $R^M$ applied to $\xi$ and $X$ is equal to  
\begin{eqnarray*}
R^M(\xi,X)\xi&=&-\nabla^M_\xi\nabla^M_X\xi+\nabla^M_X\kappa+\nabla^M_{[\xi,X]}\xi\\
&=&-\nabla^M_\xi h(X)+\nabla^M_X\kappa-g(\kappa,X)\kappa.
\end{eqnarray*}
The last equality comes from the fact that $[\xi,X]=g([\xi,X],\xi)\xi=-g(\kappa,X)\xi,$ as $X$ is foliated. Now using the O'Neill formula for Riemannian flows \cite[Eq. 4.4]{Hab} 
$$\nabla^M_\xi Y=\nabla_\xi Y+h(Y)-g(\kappa,Y)\xi,$$
for all $Y\in \Gamma(Q)$ and the fact that the tensor $h$ is a basic tensor as mentioned before, the curvature reduces to  
$$R^M(\xi,X)\xi=-h^2(X)+g(\kappa,h(X))\xi+\nabla^M_X\kappa-g(\kappa,X)\kappa.$$
This finishes the proof of the lemma. 
\hfill$\square$

\noindent At a point $x\in M,$ let us denote by $\gamma^M_0(x)$ and $\gamma^M_1(x)$ the smallest and largest eigenvalues of the symmetric tensor $R^M:\Lambda^2(Q)\rightarrow \Lambda^2(Q)$ defined by $g(R^M(X\wedge Y),Z\wedge W):=R^M_{XYZW}$ for $X,Y,Z,W\in \Gamma(Q).$ Again using the O'Neill formulas in \cite{O}, this curvature term is related to the one on the normal bundle $Q$ by the following relation: for all sections $X,Y,Z,W$ of $Q$, we have
\begin{equation} \label{eq:gh}
R^M_{XYZW}=R_{XYZW}-2g(h(X),Y)g(h(Z),W)+g(h(Y),Z)g(h(X),W)+g(h(Z),X)g(h(Y),W).
\end{equation}

\noindent Therefore according to Equation \eqref{eq:gh}, the curvature of $Q$ splits into $R_{\mbox{ext}}$ and $R_{\mbox{res}},$ where we set
$$g(R_{\mbox{ext}}(X\wedge Y), Z\wedge W)=2g(h(X),Y)g(h(Z),W)-g(h(Y),Z)g(h(X),W)-g(h(Z),X)g(h(Y),W)$$ and $$g(R_{\mbox{res}}(X\wedge Y), Z\wedge W)= R^M_{XYZW}.$$
Hence, Equation \eqref{curvature} can be written in the following way
\begin{equation}\label{eq:decomposition}
\mathcal{B}^{[p]}=\mathcal{B}^{[p]}_{\mbox{ext}}+\mathcal{B}^{[p]}_{\mbox{res}},
\end{equation}
where
\begin{equation}\label{eq:bext}
\langle \mathcal{B}^{[p]}_{\mbox{ext}} \omega , \varphi \rangle = \dfrac{1}{4} \sum_{r,s=1}^{\binom{q}{2}} \langle R_{\mbox{ext}} \psi_r , \psi_s \rangle \langle [\hat{\psi}_r, \omega],  [\hat{\psi}_s, \varphi]    \rangle 
\end{equation}
and 
$$\langle \mathcal{B}^{[p]}_{\mbox{res}} \omega , \varphi \rangle = \dfrac{1}{4} \sum_{r,s=1}^{\binom{q}{2}} \langle R_{\mbox{res}} \psi_r , \psi_s \rangle \langle [\hat{\psi}_r, \omega],  [\hat{\psi}_s, \varphi]    \rangle .$$
Choosing an orthonormal basis of eigenvectors of $R_{\mbox{res}},$ we get the pointwise estimate  
\begin{equation}\label{inequalitycurvature}
p(q-p) \gamma_0^M(x) \le \mathcal{B}^{[p]}_{\mbox{res}} \le p(q-p)\gamma_1^M(x).
\end{equation} 
Here, we use the fact that for any form $\omega\in \Lambda^p(Q)$, one has the formula 
\begin{equation}\label{norm}
\frac{1}{4} \sum_{r=1}^{\binom{q}{2}}|[\hat{\psi}_{r}, \omega]|^2=p(q-p)|\omega|^2
\end{equation}
which follows from \cite[Lem. 18]{S}. 

\noindent In order to find a lower bound of the term $\langle \mathcal{B}^{[p]}_{\mbox{ext}} \omega , \omega \rangle,$ we will compute the eigenvalues of $R_{\mbox{ext}}$ in terms of the eigenvalues of the tensor $h.$ 

{\noindent\bf Computation of the eigenvalues of the tensor $R_{\mbox{ext}}:$} Let us first check the case where $q$ is even, say $q=2m$. Since the tensor $h$ is skew-symmetric and a basic form, we can always find a local basic orthonormal frame $\{e_i\}_{i=1, \dots , q}$ of $Q$ such that the matrix of $h$ in this basis can be written as

\[
\displaystyle \left(
\begin{array}{ccccc}
  \left ( \begin{array}{cc}0&-b_1\\b_1& 0\end{array}\right ) & \raisebox{-5pt}{{\fontsize{18} {40}\selectfont \mbox{{$0$}}}} &\dots&&\raisebox{-5pt}{{\fontsize{18} {40}\selectfont \mbox{{$0$}}}} \\ [3ex]
  \raisebox{-3pt}{{\fontsize{18} {40}\selectfont \mbox{{$0$}}}} & \left ( \begin{array}{cc}0&-b_2\\b_2& 0\end{array}\right ) & &\ddots &  \vdots\\[3ex]
    \vdots &  &\ddots & & \\
     &\ddots & & &\raisebox{-3pt}{{\fontsize{18} {40}\selectfont \mbox{{$0$}}}}\\[3ex]
    \raisebox{-3pt}{{\fontsize{18} {40}\selectfont \mbox{{$0$}}}}&\dots&&\raisebox{-3pt}{{\fontsize{18} {40}\selectfont \mbox{{$0$}}}}&    \left ( \begin{array}{cc}0&-b_m\\b_m& 0\end{array}\right ) \\
\end{array}
\displaystyle  \right)
\]
where $b_1,\cdots ,b_m$ are smooth basic functions on $M$ chosen in a way such that $|b_1| \le |b_2| \le \dots \le |b_m|$. That is, $h(e_{2i-1})=b_i e_{2i}$ and $h(e_{2i})=-b_i e_{2i-1}$ for all $i=1,\cdots, m.$ Depending on the different choices of indices, we will now compute $R_{\mbox{ext}}.$ For all $i , j, k, l \in \{1, \dots, q\},$ we have
$$
\left\{
\begin{array}{lll}
\medskip
g(R_{\mbox{ext}}(e_{2i-1}\wedge e_{2i}), e_{2i-1}\wedge e_{2i})&=&3b_i^2 \\
\medskip
g(R_{\mbox{ext}}(e_{2i-1}\wedge e_{2j-1}), e_{2k-1}\wedge e_{2k})&=& 2b_ib_k \mbox{ for } k \ne i\\
\medskip
g(R_{\mbox{ext}}(e_{2i-1}\wedge e_{2j-1}), e_{2k}\wedge e_{2l})&=& -b_i b_j \delta_{jk}\delta_{il} + b_i b_j \delta_{ik}\delta_{jl}  \\
g(R_{\mbox{ext}}(e_{2i-1}\wedge e_{2j}), e_{2k-1}\wedge e_{2l})&=&2b_i b_k \delta_{ij}\delta_{kl} + b_i b_j \delta_{jk}\delta_{il}.  \\
\end{array}
\right.
$$
The other terms are all equal to zero. Therefore, in the basis $\{e_i\wedge e_j\}_{1 \le i <j \le 2m}$, arranged as follows 
$$\{e_{2i-1} \wedge
e_{2i}\}_{1 \le i \le m},\,\,\,\, \{e_{2i-1} \wedge e_{2j-1}, e_{2i}
\wedge e_{2j}\}_{1 \le i < j \le m},\,\,\,\, \{e_{2i-1} \wedge e_{2j}, e_{2i}
\wedge e_{2j-1}\}_{1 \le i < j \le m}$$

the tensor $R_{\mbox{ext}}$ is a block diagonal matrix having diagonal blocks matrices $D, D_{i,j}, -D_{i,j}$, for $1\leq i<j\leq m$ 
where:
\begin{itemize}
\item $D$ is the matrix representation of the restriction of
$R_{\mbox{ext}}$ to the subspace generated by $\{e_{2i-1} \wedge
e_{2i}\}_{1 \le i \le m}$ and is given by
$$
D=\begin{pmatrix}
  3b_1^2 & 2b_1 b_2 & \hdots & 2b_1b_m \\
  2b_1 b_2 & 3b_2^2 & \hdots & 2b_2b_m \\
 \vdots &  & \ddots &  \\
  2b_1b_m & 2b_2b_m & \hdots & 3b_m^2 \\
\end{pmatrix}.
$$
\item $D_{i,j}$ is the
matrix representation of the restriction of $R_{\mbox{ext}}$ to
the subspace generated by $\{e_{2i-1} \wedge e_{2j-1}, e_{2i}
\wedge e_{2j}\}$ which is given by
$$
\begin{pmatrix}
  0 & b_ib_j \\
  b_i b_j & 0 \\
\end{pmatrix}.
$$
\item The last block $-D_{i,j}$ is the matrix representation of the restriction of $R_{\mbox{ext}}$
to the subspace generated by $\{e_{2i-1} \wedge e_{2j}, e_{2i}\wedge e_{2j-1}\}$.
\end{itemize}
\noindent We notice that by a straightforward computation one can prove that the choice of the basis does not change the orientation of the normal bundle.

\noindent One can easily check that the eigenvalues of the matrices $D_{i,j}$ are $\pm b_{i}b_{j}$ with unit
eigenvectors $\theta_{ij}^{\pm}= \frac{1}{\sqrt{2}}(e_{2i-1} \wedge
e_{2j-1}\pm e_{2i} \wedge e_{2j}).$ Also the eigenvalues of the matrices $-D_{i,j}$ are $\pm b_{i}b_{j}$ with unit
eigenvectors given by $\rho_{ij}^\mp= \frac{1}{\sqrt{2}}(e_{2i-1} \wedge
e_{2j}\mp e_{2i} \wedge e_{2j-1}).$ The eigenvalues of the matrix $D$ are not easy to compute but we know that they are all nonnegative since $\langle DX, X \rangle =  \sum_{i=1}^m b_i^2 X_i^2 +2(\sum_{i=1}^m b_i X_i)^2\geq 0$ for any vector $X.$

\noindent In conclusion, the eigenvalues $\{\lambda_r\}_{r=1,\cdots,\binom{q}{2}}$ of the tensor $R_{\mbox{ext}}$ consist of three families ($q$ is even): 

\begin{itemize} 
\item {\bf Type I :} The eigenvalues are $\pm b_{i}b_{j}$ ($i<j$) with unit eigenvectors $\theta_{ij}^{\pm}= \frac{-1}{\sqrt{2}}(e_{2i-1} \wedge
e_{2j-1}\pm e_{2i} \wedge e_{2j})$ 
\item {\bf Type II :} The eigenvalues are $\pm b_{i}b_{j}$ ($i<j$) with unit eigenvectors given by $\rho_{ij}^\mp= \frac{1}{\sqrt{2}}(e_{2i-1} \wedge
e_{2j}\mp e_{2i} \wedge e_{2j-1}).$ 
\item {\bf Type III :} The eigenvalues are those of the matrix $D$ which are all nonnegative and the eigenvectors are in the subspace generated by $\{e_{2i-1}\wedge e_{2i}\}_{i=1\cdots,m}.$ 
\end{itemize}

\noindent The case where $q$ is odd can be treated in a similar way as the even case but an additional direction $e_0$ is involved corresponding to the eigenvalue $0$ of $h.$ Since $g(R_{\mbox{ext}}(e_0\wedge X),Y\wedge Z)=0$ for every $X,Y,Z\in \Gamma(Q),$ we deduce that the eigenvalues of $R_{\mbox{ext}}$ consist of families of type {\bf I, II, III} (the same as defined above) and {\bf IV}, where in the last family $0$ is an eigenvalue and the corresponding eigenvector is in the subspace generated by $\{e_0\wedge e_i\}_{i=1,\cdots, 2m}.$

{\bf Lower bound of the term $\langle \mathcal{B}^{[p]}_{\mbox{ext}} \omega , \omega\rangle$:} Let us denote by $\tilde{\lambda}_r$ ($1 \le r \le m$) the eigenvalues of the matrix $D$ and let $\{\tilde{\theta}_r\}$ be an orthonormal family of
eigenvectors associated with the eigenvalues $\tilde{\lambda}_r$. Then we have the estimate,

\begin{eqnarray}\label{inequalitycurvatureext}
\langle \mathcal{B}^{[p]}_{\mbox{ext}} \omega , \omega\rangle & =&
\dfrac{1}{4} \sum_{1\le i <j \le m} b_i b_j( |[{\theta}_{ij}^+,
\omega]|^2  +|[{\rho}_{ij}^-, \omega]|^2 ) \nonumber- \dfrac{1}{4} \sum_{1\le i <j \le m} b_i b_j( |[{\theta}_{ij}^-,
\omega]|^2 +|[{\rho}_{ij}^+, \omega]|^2 )  \nonumber \\  &&+
\dfrac{1}{4} \sum_{r=1}^m \tilde{\lambda}_r |[\tilde{\theta}_{r},
\omega]|^2\nonumber\\
&\ge & \dfrac{1}{4} \sum_{1\le i <j \le m} b_i b_j( |[{\theta}_{ij}^+,
\omega]|^2  +|[{\rho}_{ij}^-, \omega]|^2 )    - \dfrac{1}{4}
\sum_{1\le i <j \le m} b_i b_j( |[{\theta}_{ij}^-, \omega]|^2
+|[{\rho}_{ij}^+, \omega]|^2 )   \nonumber\\
& \ge & -\dfrac{1}{4} b_m^2 \left( \sum_{1\le i <j \le m} \left(
|[{\theta}_{ij}^+, \omega]|^2  +|[{\rho}_{ij}^-, \omega]|^2+|[{\theta}_{ij}^-, \omega]|^2
+|[{\rho}_{ij}^+, \omega]|^2 \right) \right)  \nonumber\\
& \bui{\ge}{\eqref{norm}} & -p(q-p)b_m^2 |\omega|^2.
\end{eqnarray}

\noindent Hence, we arrive at the following result: 

\begin{thm}\label{curvatureestimate}
Let $(M^{2m+1},g)$ be a Riemannian manifold endowed with a
Riemannian flow given by a unit vector field $\xi$ of
codimension $q$. For any number $p$ such that $1\leq p \leq q-1$ and a basic $p$-form $\omega,$ we have 
\begin{equation*} 
\langle \mathcal{B}_{\rm{ext}}^{[p]}\omega,\omega\rangle \geq -p(q-p)b_m^2|\omega|^2,
\end{equation*}
with $m=[\frac{q}{2}].$ If the equality is attained for some $p\in \{1,\cdots, m\}$, then $|b_1|=\cdots=|b_m|.$ If $m=1$ and the equality is attained, then $b_1=0.$ 
\end{thm} 

\noindent Before proving the theorem, let us give some direct consequences.

\begin{cor}\label{curvatureestimate3}
Let $(M,g)$ be a Riemannian manifold endowed with a
Riemannian flow given by a unit vector field $\xi$ of
codimension $q$. For any number $p$ such that $1\leq p \leq q-1$ and any basic $p$-form $\omega,$ we have
\begin{equation*} 
\langle \mathcal{B}^{[p]}\omega,\omega\rangle \geq p(q-p)(\gamma_0^M-b_m^2)|\omega|^2,
\end{equation*}
where $\gamma_0^M$ is the lowest eigenvalue of the curvature operator of $M$ restricted to $Q$ and $m=[\frac{q}{2}].$ If $m>1$ and the equality is attained for some $p\in \{1,\cdots, m\},$ then $|b_1|=\cdots=|b_m|.$ If $m=1$ and the equality is attained, then $b_1=0.$   
\end{cor}

\noindent The inequality in Corollary \ref{curvatureestimate3} is obtained by adding the estimate in Theorem \ref{curvatureestimate} to the l.h.s. of Inequality \eqref{inequalitycurvature}. One can easily check that for the Hopf fibration $\mathbb{S}^{2m+1}\rightarrow \mathbb{C}{\rm P}^m$ for $m>1,$ the K\"ahler form $\Omega$ on $\mathbb{C}{\rm P}^m,$ which is a parallel basic $2$-form, satisfies the equality of the above theorem (here $\gamma_0^M=b_m^2=1$). Also on the Riemmanian product $\mathbb{S}^1\times \mathbb{S}^{2m+1}$ for $m>1,$ when one considers the flow defined by the unit vector field $\xi:=\frac{1}{\sqrt{2}}(\xi_1+\xi_2)$ where $\xi_1$ is the unit parallel vector field on $\mathbb{S}^1$ and $\xi_2$ is the unit Killing vector field that defines the Hopf fibration, the K\"ahler form on $\mathbb{C}{\rm P}^m$ is transversally parallel. In this case, the equality is attained since $\gamma_0^M=b_m^2=\frac{1}{2}.$ We point out that the converse of Corollary \ref{curvatureestimate3} does not hold in general. Indeed, consider the Riemannian fibration $\mathbb{S}^1\times \mathbb{S}^{2m+1}\rightarrow \mathbb{S}^1\times \mathbb{C}{\rm P}^m$ and let $\Omega$ be again the K\"ahler form on $\mathbb{C}{\rm P}^m.$ Here $|b_1|=\cdots=|b_m|=1$ and $\gamma_0^M=0$ which gives the strict inequality.

\noindent When the term in the lower bound of Corollary \ref{curvatureestimate3} is positive, we get the following rigidity result: 

\begin{cor} \label{corrigidity}
Let $(M,g)$ be a compact Riemannian manifold endowed with a Riemannian flow given by a unit vector field $\xi$ of
codimension $q$. If $\gamma_0^M\geq b_m^2$ and $\kappa$ is basic-harmonic, then every harmonic basic $p$-form is transversally parallel. If the strict inequality holds, then $H_b^s(\mathcal{F})=\{0\}$ for any $s\in \{1,\cdots,q-1\}.$
\end{cor}

\noindent The proof of this corollary uses the first statement of Proposition \ref{estimate}. Another direct consequence of Corollary \ref{curvatureestimate3} that characterizes minimal Riemannian flow on round spheres is the following (see \cite{GG})

\begin{cor}\label{corrigidity2}
Let $\mathbb{S}^n$ be the round sphere of constant sectionnal curvature $1$ and assume that it is endowed with a minimal Riemannian flow. Then, the O'Neill tensor is transversally parallel and the flow defines a Sasakian structure on $\mathbb{S}^n.$ 
\end{cor} 

{\noindent \it Proof of Corollary \ref{corrigidity2}:} As the curvature on the sphere $\mathbb{S}^n$ is given for all vector fields $X,Y,Z$ by $R^M(X,Y)Z=g(X,Z)Y-g(Y,Z)X,$ one deduces directly from Lemma \ref{computationcurvature} that $h^2(X)=-X$ for all $X\in \Gamma(Q),$ that is $|b_1|=\cdots=|b_m|=1$. In the same way, using the fact that for $X,Y,Z\in \Gamma(Q),$ we have \cite{O}
$$g(R^M(X,Y)\xi,Z)=g(-(\nabla_X h)Y+(\nabla_Y h)X,Z),$$
one can also get that $(\nabla_X h)Y=(\nabla_Y h)X.$ Recall here that $\nabla$ is the transversal Levi-Civita connection extended to forms. Therefore, the divergence of $h$ (with respect to the normal bundle) vanishes since
$$ (\delta h)(X)=-\sum_{i=1}^{n-1}(\nabla_{e_i}h)(e_i,X)=\sum_{i=1}^{n-1}(\nabla_{e_i}h)(X,e_i)=\sum_{i=1}^{n-1}(\nabla_{X}h)(e_i,e_i)=0.$$
Hence, the basic $2$-form $\Omega:=-\frac{1}{2}d\xi=g(h\cdot,\cdot)$ is closed and coclosed and thus a basic-harmonic. Now, Corollary \ref{corrigidity} allows to deduce that it is transversally parallel. This ends the proof. 
\hfill$\square$




\noindent Using the second statement in Proposition \ref{estimate}, one can deduce the following estimate

\begin{thm}\label{estimateeigenvalues}
Let $(M,g)$ be a compact Riemannian manifold endowed with a minimal
Riemannian flow given by a unit vector field $\xi$ of
codimension $q$. Let $p$ be any integer number such that $1\leq p\leq m$ with $m=[\frac{q}{2}].$ Then the first eigenvalue of the basic Laplacian acting on basic $p$-forms satisfies 
$$\lambda_{1,p} \ge p(q-p+1)(\gamma_M +\beta_M^1),$$
where $\gamma_M={\rm inf}_M(\gamma_0^M)$ is a lower bound of the curvature operator on $M$ restricted to $Q$ and $\beta_M^1={\rm inf}_M(-b_m^2)$ is the lowest eigenvalue of the symmetric tensor $h^2.$ If moreover the equality is attained, then $M$ is isometric to the quotient of $\mathbb{R}\times \Sigma$ by some fixed-point-free cocompact discrete subgroup $\Gamma\subset \mathbb{R}\times {\rm SO}_{q+1},$ where $\Sigma$ is a compact simply connected manifold of positive curvature. 
\end{thm} 

{\noindent \bf Remarks.} 
\begin{enumerate} 
\item In the equality case of the estimate in Theorem \ref{estimateeigenvalues}, the O'Neill tensor vanishes. Therefore, the basic Laplacian on $M$ restricts to the usual Laplacian on the manifold $\Sigma$ and thus the first eigenvalue on $\Sigma$ satisfies the equality case in the Gallot-Meyer estimate \cite[Thm. 6.13]{GM}. In view of the remark after Theorem \ref{estimate} and if $p$ is chosen such that $p<\frac{q}{2},$ we deduce that $d\omega=0$ where $\omega$ is an eigenform associated with the first eigenvalue. If $p=2$ and $q>4,$ the form $\alpha=\delta\omega$ is a coclosed $1$-form which is still an eigenform of the Laplacian (the form $\alpha$ does not vanish since this would imply that $\omega$ vanishes). Hence, by a result of S. Tachibana \cite[Thm. 3.3]{Ta} the manifold $\Sigma$ is either isometric to a Sasakian manifold or to a round sphere with constant curvature.  
\item By the result in \cite{BW}, the manifold $\Sigma$ is a spherical space form. In case $\Sigma$ is isometric to a round sphere, the group $\Gamma=\pi_1(M)$ preserves the orthogonal splitting $T_{(t,x)}\widetilde M=\mathbb{R}\oplus T_x \mathbb{S}^q$ (the vertical distribution $\mathbb{R}$ is the kernel of the Ricci tensor), as it is acting by isometries on the universal cover $\widetilde M.$ Therefore the fundamental group is embedded in the product ${\rm Isom}_+(\mathbb{R})\times {\rm Isom}_+(\mathbb{S}^q)$ where ${\rm Isom}_+$ is the group of isometries that preserve the orientation of the corresponding manifold. For $q$ even, we deduce that $\Gamma\simeq \mathbb{Z}$ and that it acts as $(t,x)\rightarrow (t+a, A(x))$ for some $(a,A)\in \mathbb{R}^*\times {\rm SO}(q+1).$ For $q$ odd, the group $\Gamma$ is not necessarily isomorphic to $\mathbb{Z},$ since one might consider the group $\Gamma=\mathbb{Z}\times \Gamma_2$ where $\Gamma_2$ is a finite subgroup of ${\rm SO}(q+1)$ consisting of rotations in orthogonal $2$-planes in $\mathbb{R}^{q+1}.$  
\end{enumerate}

\noindent Let us now proceed with the proofs of the equality case of Theorems \ref{curvatureestimate} and \ref{estimateeigenvalues}.

{\noindent \it Proof of Theorem \ref{curvatureestimate}:} First, we discuss the case where $q=2m>2.$ If the equality is attained in \eqref{inequalitycurvatureext}, then two cases may occur: Either for all $(i,j)$ one of the Lie bracket coefficients of $b_i b_j$ in the first line of \eqref{inequalitycurvatureext} does not vanish and in this case we get $|b_1|=\cdots=|b_m|$ or there exist $i$ and $j$ with $i<j$ and such that all the coefficients vanish, that is 
\begin{equation}\label{equalitybracket}
[{\theta}_{ij}^\pm,\omega]=[{\rho}_{ij}^\pm, \omega]=0.
\end{equation}
 
\noindent Let us check that the second case gives also the statement of the theorem. First, we get a description of the form $\omega$ that we put it in the following lemma: 

\begin{lemma} \label{description} Assume that there exist $i,j$ such that Equalities \eqref{equalitybracket} hold. Then, there exist basic forms $\omega_1$ and $\omega_2$ such that 
$$\omega=e_{2i-1}\wedge e_{2i}\wedge e_{2j-1}\wedge e_{2j}\wedge \omega_1+\omega_2,$$
with 
$$
\left\{\begin{array}{ll}
e_{2i-1}\lrcorner\, \omega_1=e_{2i}\lrcorner\, \omega_1=0\\
e_{2j-1}\lrcorner\, \omega_1=e_{2j}\lrcorner\, \omega_1=0,
\end{array}\right.
$$
The same equalities hold for $\omega_2$.
\end{lemma}

\noindent{\it Proof.} By adding (and substracting) the brackets $[{\theta}_{ij}^+,\omega]$ and $[{\theta}_{ij}^-,\omega]$ together, as well as $[{\rho}_{ij}^+, \omega]$ and $[{\rho}_{ij}^-, \omega],$ we deduce the following equations 
$$[e_{2i-1}\wedge e_{2j-1},\omega]=[e_{2i}\wedge e_{2j},\omega]=[e_{2i-1}\wedge e_{2j},\omega]=[e_{2i}\wedge e_{2j-1},\omega]=0.$$
Now, using Lemma \ref{bracketcomputation1} for each of the above brackets, the previous equations reduce to the following system 

$$
\left\{\begin{array}{ll}
e_{2j-1}\wedge (e_{2i-1}\lrcorner \omega)=e_{2i-1}\wedge (e_{2j-1}\lrcorner \omega)\\\\
e_{2j}\wedge (e_{2i}\lrcorner \omega)=e_{2i}\wedge (e_{2j}\lrcorner \omega)\\\\
e_{2j}\wedge (e_{2i-1}\lrcorner \omega)=e_{2i-1}\wedge (e_{2j}\lrcorner \omega)\\\\
e_{2j-1}\wedge (e_{2i}\lrcorner \omega)=e_{2i}\wedge (e_{2j-1}\lrcorner \omega).
\end{array}\right.
$$
In order to solve this system, we take the interior product of the first equation with $e_{2i-1}$ (resp. with $e_{2j-1}$) to get that 
$$e_{2i-1}\lrcorner \omega=e_{2j-1}\wedge \beta_0  \quad\text{and}\quad e_{2j-1}\lrcorner \omega=e_{2i-1}\wedge \beta_1,$$ 
where $\beta_0$ (resp. $\beta_1)$ is a form that does not contain neither $e_{2i-1}$ nor $e_{2j-1}$. The same can be done for the third equation with respect to $e_{2i-1}$ and $e_{2j}$ to obtain 
$$e_{2i-1}\lrcorner \omega=e_{2j}\wedge \beta_3  \quad\text{and}\quad e_{2j}\lrcorner \omega=e_{2i-1}\wedge \beta_4,$$
for some $\beta_3,\beta_4.$ Comparing the above equations and using the fact that the general solution of an equation of type $X\wedge \alpha=Y\wedge \beta$ where $X$ and $Y$ are orthogonal and $X\lrcorner \alpha=Y\lrcorner \beta=0$ is given by $\alpha=Y\wedge (X\lrcorner \beta),$ we conclude that $\beta_0$ should be of the form $e_{2j}\wedge \beta_5$ for some form $\beta_5.$ The same technique can be used for the second and forth equations in the system. This allows to finish the proof of the lemma by using the fact that the general solution of an equation of the form $X\lrcorner \omega=\alpha$ is $\omega=X\wedge \alpha+\beta$ where $X\lrcorner \beta=0.$
\hfill$\square$   

\noindent We now proceed with the proof of Theorem \ref{curvatureestimate}. According to Lemmas \ref{description}, \ref{bracketcomputation} and to Equality \eqref{eq:bext}, we set $\Phi:=e_{2i}\wedge e_{2j-1}\wedge e_{2j}\wedge \omega_1$ and we write
 
\begin{eqnarray}\label{eq:inequality ext} 
\langle\mathcal{B}_{{\rm ext}}^{[p]}\omega,\omega\rangle&=&\langle\mathcal{B}_{{\rm ext}}^{[p]}(e_{2i-1}\wedge\Phi),e_{2i-1}\wedge\Phi\rangle+2\langle\mathcal{B}_{{\rm ext}}^{[p]}(e_{2i-1}\wedge\Phi),\omega_2\rangle+\langle\mathcal{B}_{{\rm ext}}^{[p]}\omega_2,\omega_2\rangle\nonumber\\
&=&\frac{1}{4}\sum_{r=1}^{\binom{q}{2}} \lambda_r|[\hat\theta_r,e_{2i-1}\wedge\Phi]|^2+\frac{1}{2}\sum_{r=1}^{\binom{q}{2}} \lambda_r \langle[\hat\theta_r,e_{2i-1}\wedge\Phi],[\hat\theta_r,\omega_2]\rangle+\langle\mathcal{B}_{{\rm ext}}^{[p]}\omega_2,\omega_2\rangle\nonumber\\
&=&\frac{1}{4}\sum_{r=1}^{\binom{q}{2}} \lambda_r|[\hat\theta_r,\Phi]|^2+\sum_{r=1}^{\binom{q}{2}} \lambda_r|e_{2i-1}\lrcorner \hat\theta_r|^2|\Phi|^2+\sum_{r=1}^{\binom{q}{2}} \lambda_r\langle e_{2i-1}\cdot[\hat\theta_r,\Phi],(e_{2i-1}\lrcorner \hat\theta_r)\cdot \Phi\rangle\nonumber\\
&&+\frac{1}{2}\sum_{r=1}^{\binom{q}{2}} \lambda_r \langle e_{2i-1}\cdot[\hat\theta_r,\Phi],[\hat\theta_r,\omega_2]\rangle+\sum_{r=1}^{\binom{q}{2}} \lambda_r \langle (e_{2i-1}\lrcorner\hat\theta_r)\cdot \Phi,[\hat\theta_r,\omega_2]\rangle+\langle\mathcal{B}_{{\rm ext}}^{[p]}\omega_2,\omega_2\rangle\nonumber\\
&=&\langle\mathcal{B}_{{\rm ext}}^{[p-1]}\Phi,\Phi\rangle+\sum_{r=1}^{\binom{q}{2}} \lambda_r|e_{2i-1}\lrcorner \hat\theta_r|^2|\Phi|^2-\sum_{r=1}^{\binom{q}{2}} \lambda_r\langle [\hat\theta_r,\Phi],e_{2i-1}\cdot (e_{2i-1}\lrcorner \hat\theta_r)\cdot \Phi\rangle\nonumber\\
&&+\frac{1}{2}\sum_{r=1}^{\binom{q}{2}} \lambda_r \langle e_{2i-1}\wedge[\hat\theta_r,\Phi],[\hat\theta_r,\omega_2]\rangle+\sum_{r=1}^{\binom{q}{2}} \lambda_r \langle (e_{2i-1}\lrcorner\hat\theta_r)\wedge \Phi,[\hat\theta_r,\omega_2]\rangle+\langle\mathcal{B}_{{\rm ext}}^{[p]}\omega_2,\omega_2\rangle.\nonumber\\
\end{eqnarray}

\noindent Here, we recall that $\{\lambda_r\}$ are the eigenvalues of the tensor $R_{\mbox{ext}}$ and $\{\hat\theta_r\}$ are the corresponding dual eigenvectors found previously. In the following, we will compute each sum separately with respect to each family of eigenvalues of type {\bf (I), (II)} and {\bf (III)} that we already find. For this, we denote by ${\bf S}_1,{\bf S}_2,{\bf S}_3$ and ${\bf S}_4$ the respective sums.

\noindent{\bf Type I :} In the following, we shall prove that ${\bf S}_1,{\bf S}_2,{\bf S}_3$ and ${\bf S}_4$ all vanish with respect to an orthonormal basis of type {\bf I}. In fact, as we have that 
\begin{equation}\label{eq:interiorproduct}
e_{s}\lrcorner \theta_{kl}^{\pm}=\frac{-1}{\sqrt{2}}(\delta_{s2k-1} e_{2l-1}-\delta_{s2l-1} e_{2k-1}\pm\delta_{s2k}e_{2l}\mp\delta_{s2l}e_{2k}),
\end{equation} 
we first deduce that $|e_{2i-1}\lrcorner \theta_{kl}^{\pm}|^2=\frac{1}{2}$ if $i=k$ or $i=l$ and thus ${\bf S}_1$ is zero (the sum of all the eigenvalues). Second, from Lemma \ref{bracketcomputation1}, we have that 
\begin{equation}\label{eq:liebracket}
[\theta_{kl}^{\pm},\Theta]=\frac{-2}{\sqrt{2}}\left(e_{2l-1}\wedge (e_{2k-1}\lrcorner\Theta)-e_{2k-1}\wedge (e_{2l-1}\lrcorner\Theta)\pm e_{2l}\wedge (e_{2k}\lrcorner\Theta)\mp e_{2k}\wedge (e_{2l}\lrcorner\Theta)\right),
\end{equation}
for any form $\Theta.$ Therefore, we get that 
$$(e_{2i-1}\lrcorner \theta_{kl}^{\pm})\lrcorner [\theta_{kl}^{\pm},\omega_2]=
\left\{\begin{array}{ll}
\pm e_{2i}\wedge (e_{2k-1}\lrcorner e_{2k}\lrcorner\omega_2) &\quad{\rm for}\quad i=l\\\\
\pm e_{2i}\wedge (e_{2l-1}\lrcorner e_{2l}\lrcorner\omega_2) &\quad{\rm for}\quad i=k
\end{array}\right.$$
(up to a factor $\frac{-1}{\sqrt{2}}$) which, by taking the scalar product with $\Phi,$ gives that ${\bf S}_4=0.$ Here we used the fact that $\omega_2$ does not contain any factor in $e_i$ and $e_j.$ For the sum ${\bf S}_3$, we first compute 
$$e_{2i-1}\lrcorner [\theta_{kl}^{\pm},\omega_2]=\frac{-2}{\sqrt{2}}\left(\delta_{il} e_{2k-1}\lrcorner\omega_2-\delta_{ik} e_{2l-1}\lrcorner\omega_2\right).$$ 
Hence, the term (up to the factor $\frac{-2}{\sqrt{2}}$) 
$$
\langle [\theta_{kl}^{\pm},\Phi],e_{2i-1}\lrcorner [\theta_{kl}^{\pm},\omega_2]\rangle=\left\{\begin{array}{ll}
\langle [\theta_{ki}^{\pm},\Phi],e_{2k-1}\lrcorner \omega_2\rangle &\quad{\rm for}\quad i=l\\\\
-\langle [\theta_{il}^{\pm},\Phi],e_{2l-1}\lrcorner \omega_2\rangle &\quad{\rm for}\quad i=k
\end{array}\right.$$ 
also vanishes by Equation \eqref{eq:liebracket} (replace $\Theta$ by $\Phi$ and $l$ or $k$ by $i$).  Hence ${\bf S}_3=0.$ Now, we are left with the sum ${\bf S}_2$ that we shall prove that it vanishes as well. Indeed, we write   
\begin{eqnarray*} 
{\bf S}_2&=&\sum_{k<l}b_kb_l\langle [\theta_{kl}^+,\Phi],e_{2i-1}\cdot (e_{2i-1}\lrcorner \theta_{kl}^+)\cdot \Phi\rangle-\sum_{k<l}b_kb_l\langle [\theta_{kl}^-,\Phi],e_{2i-1}\cdot (e_{2i-1}\lrcorner \theta_{kl}^-)\cdot \Phi\rangle\\
&\bui{=}{\eqref{eq:interiorproduct}}& \sum_{i<l}b_ib_l\langle [\theta_{il}^+-\theta_{il}^-,\Phi],e_{2i-1}\cdot e_{2l-1}\cdot \Phi\rangle-\sum_{k<i}b_kb_i\langle [\theta_{ki}^+-\theta_{ki}^-,\Phi],e_{2i-1}\cdot e_{2k-1}\cdot \Phi\rangle.
\end{eqnarray*}
Now from the expression of the vector fields $\theta_{kl}^+$ and $\theta_{kl}^-$ and using again Lemma \ref{bracketcomputation}, we have that 
\begin{eqnarray*} 
\langle[\theta_{il}^+-\theta_{il}^-,\Phi],e_{2i-1}\cdot e_{2l-1}\cdot \Phi\rangle &=&\frac{-2}{\sqrt{2}}\langle[e_{2i}\wedge e_{2l},\Phi],e_{2i-1}\cdot e_{2l-1}\cdot \Phi\rangle\\
&=&\frac{-4}{\sqrt{2}}\langle e_{2l}\wedge (e_{2i}\lrcorner\Phi)-e_{2i}\wedge (e_{2l}\lrcorner\Phi),e_{2i-1}\cdot e_{2l-1}\cdot \Phi\rangle\\
&=&\frac{-4}{\sqrt{2}}\langle e_{2l}\wedge (e_{2i}\lrcorner\Phi),e_{2l-1}\lrcorner (e_{2i-1}\wedge \Phi)\rangle\\
&=&\frac{-4}{\sqrt{2}}\langle e_{2l-1}\wedge e_{2l}\wedge (e_{2i}\lrcorner\Phi),e_{2i-1}\wedge \Phi\rangle=0,
\end{eqnarray*}
which means that the first sum vanishes. By interchanging the roles of $i$ and $l,$ we also deduce that the second sum ${\bf S}_2$ is zero. 

\noindent{\bf Type II :} The computation can be done in the same way as for type {\bf I} and shows that all of the sums vanish.

\noindent{\bf Type III :} Recall that in this case, the eigenvectors of $R_{{\rm ext}}$ are in the subspace generated by $\{e_{2k-1}\wedge e_{2k}\}_{k=1\cdots,m}.$ Hence any eigenvector $\tilde\theta_r$ ($1\leq r\leq m$) can be written as $\tilde\theta_r=\sum_{k=1}^m \alpha_r^k e_{2k-1}\wedge e_{2k}$ for some functions $\alpha_r^k.$ Thus, we have 
\begin{equation} \label{eq:type3}
e_{2i-1}\lrcorner \tilde\theta_r=\alpha_r^i e_{2i}. 
\end{equation}
The first sum ${\bf S}_1$ is then equal to $\sum_{r=1}^m \tilde\lambda_r(\alpha_r^i)^2|\Phi|^2,$ where $\tilde\lambda_r$ are the eigenvalues of the matrix $D$ defined before. Next, we shall prove that ${\bf S}_3$ and ${\bf S}_4$ are equal to zero. Indeed, using \eqref{eq:type3}, one can easily see that $(e_{2i-1}\lrcorner \tilde\theta_r)\wedge \Phi=0$ which gives that ${\bf S}_4=0.$ Now using Lemma \ref{bracketcomputation}, one has
$$[\tilde\theta_r,\Theta]=\sum_{k=1}^m \alpha_r^k [e_{2k-1}\wedge e_{2k},\Theta]=2 \sum_{k=1}^m \alpha_r^k (e_{2k}\wedge (e_{2k-1}\lrcorner \Theta)-e_{2k-1}\wedge (e_{2k}\lrcorner \Theta)),$$ 
for any form $\Theta.$ This gives that $e_{2i-1}\lrcorner [\tilde\theta_r,\omega_2]=0$ and thus ${\bf S}_3=0.$ Here, we used the fact that $\omega_2$ does not contain any factor in $e_i.$ The term ${\bf S}_2$ is now equal to 
\begin{eqnarray*} 
{\bf S}_2&=&\sum_{r=1}^m\lambda_r \alpha_r^i \langle [\tilde\theta_r,\Phi],e_{2i-1}\cdot e_{2i}\cdot\Phi\rangle\\
&=&2 \sum_{k,r=1}^m \lambda_r \alpha_r^i\alpha_r^k \langle(e_{2k}\wedge (e_{2k-1}\lrcorner \Phi)-e_{2k-1}\wedge (e_{2k}\lrcorner \Phi),e_{2i-1}\wedge (e_{2i}\lrcorner\Phi)\rangle\\
&=&-2 \sum_{k,r=1}^m \lambda_r \alpha_r^i\alpha_r^k \langle e_{2k-1}\wedge (e_{2k}\lrcorner \Phi),e_{2i-1}\wedge (e_{2i}\lrcorner\Phi)\rangle\\
&=&-2 \sum_{k,r=1}^m \lambda_r \alpha_r^i\alpha_r^k\delta_{ik}|e_{2i}\lrcorner\Phi|^2=-2\sum_{r=1}^m \lambda_r(\alpha_r^i)^2|\Phi|^2. 
\end{eqnarray*}
Now replacing all the computations above in Equation \eqref{eq:inequality ext}, we deduce that  
\begin{eqnarray*}
-p(q-p)b_m^2|\omega|^2=\langle\mathcal{B}_{{\rm ext}}^{[p]}\omega,\omega\rangle&=&\langle\mathcal{B}_{{\rm ext}}^{[p-1]}\Phi,\Phi\rangle+3\sum_{r=1}^m \lambda_r(\alpha_r^i)^2|\Phi|^2+\langle\mathcal{B}_{{\rm ext}}^{[p]}\omega_2,\omega_2\rangle\\
&\bui{\geq}{\eqref{inequalitycurvatureext}} & -(p-1)(q-p+1)b_m^2|\Phi|^2-p(q-p)b_m^2|\omega_2|^2.
\end{eqnarray*}
Here, we use the fact that all the eigenvalues $\lambda_r$ are nonnegative. As $|\omega|^2=|\Phi|^2+|\omega_2|^2,$ the last inequality implies that either $b_m=0$ or that $\Phi=0.$  Recall here that the integer $p$ is chosen such that $1\leq p\leq m.$ The fact that the $b_i$'s are chosen in a way that $|b_1|\leq \cdots \leq |b_m|,$ then $b_m=0$ implies the statement of Theorem \ref{curvatureestimate}. We are now left with the case when $\Phi=0,$ which means by Lemma \ref{description} that $\omega=\omega_2$ with $e_{2i-1}\lrcorner \omega=e_{2i}\lrcorner \omega=e_{2j-1}\lrcorner \omega=e_{2j}\lrcorner \omega=0.$ But recall that $i$ and $j$ are chosen in a way that all the Lie bracket coefficients of $b_ib_j$ in Equation \eqref{inequalitycurvatureext} are equal to zero. Therefore the same choice holds for $i=1$ and $1\leq j\leq m,$ since otherwise we would get $|b_1|=\cdots=|b_m|.$ Hence by varying $j,$ we arrive at $X\lrcorner \omega=0$ for any $X,$ which leads to $\omega=0;$ that is a contradiction. This finishes the proof for $m>1.$

\noindent  Now, we discuss the equality when $q$ is odd, say $q=2m+1.$ In this case, we have $[e_0\wedge e_l,\omega]=0$ for all $l=1,\cdots,2m.$ Recall here that $e_0$ is the eigenvector of $h$ that corresponds to the eigenvalue $0.$ As in the even case, either for all $(i,j)$ one of the Lie bracket coefficients of $b_i b_j$ in \eqref{inequalitycurvatureext} does not vanish and we get $|b_1|=\cdots=|b_m|$ or there exist $i$ and $j$ with $i<j$ and such that all the coefficients vanish. In the second case, Equations \eqref{equalitybracket} still hold and we get the same description as in Lemma \ref{description}. That means, we write $\omega=e_{2i-1}\wedge e_{2i}\wedge e_{2j-1}\wedge e_{2j}\wedge \omega_1+\omega_2.$ From the one hand, we take $l=2i-1$ in the equation $[e_0\wedge e_l,\omega]=0$ and make the interior product of this last identity with $e_{2i-1}$ to get after using Lemma \ref{bracketcomputation1}
 
\begin{equation}\label{odd1}
e_0\lrcorner \omega_2=0 \quad\text{and}\quad e_0\wedge \omega_1=0.
\end{equation}
From the other hand, we take $l\notin\{2i-1,2i,2j-1,2j\}$ and make the interior product of the same equation with $e_{2i-1}\wedge e_{2i}\wedge e_{2j-1}\wedge e_{2j}$ to find that 
\begin{equation}\label{odd2}
e_l\wedge (e_0\lrcorner\omega_1)=0 \quad\text{and}\quad e_0\wedge (e_l\lrcorner\omega_2)=0. 
\end{equation}

\noindent Now, the interior product of the first equation in \eqref{odd1} with $e_l$ and the second equation in \eqref{odd2} allow to deduce that $\omega_2=0.$ Therefore, we deduce that $\omega=e_{2i-1}\wedge \Phi.$ The rest of the proof carries on the same way as in the even case. We notice that the family {\bf IV} of eigenvalues does not contribute to Equation \eqref{eq:inequality ext}, since in this case all the eigenvalues are equal to zero.  

\noindent We are now left with the case when $m=1.$ As from the first line of Equation \eqref{inequalitycurvatureext} the term $\langle \mathcal{B}_{{\rm ext}}^{[1]}\omega,\omega\rangle$ is nonnegative, we then deduce that the equality in Theorem \ref{curvatureestimate} is attained if $b_1=0.$ This ends the proof. 
\hfill$\square$

{\noindent \it Proof of Theorem \ref{estimateeigenvalues}:} Assume that the estimate is realized, then the inequality in Corollary \ref{curvatureestimate3} is also attained and therefore $|b_1|=\cdots=|b_m|={\rm cst}$ for $m>1$ and $b_1=0$ for $m=1.$ In the following, we will prove that the constant should also be zero. Indeed, as $\lambda_{1,p}=p(q-p+1)(\gamma_M-{\rm cst})>0$ we deduce that $\gamma_M>{\rm cst}>0.$ Therefore from Corollary \ref{corrigidity}, we get that $H^2_b(\mathcal{F})=0.$ On the other hand, using Lemma \ref{computationcurvature}, the Ricci curvature on $M$ is equal to
$${\rm Ric}^M(\xi,\xi)=\sum_{i=1}^qR^M(\xi,e_i,\xi,e_i)=-\sum_{i=1}^q g(h^2e_i,e_i)=|h|^2=2m{\rm cst}>0,$$ 
and 
$${\rm Ric}^M(X,X)=\sum_{i=1}^qR^M(X,e_i,X,e_i)+R^M(X,\xi,X,\xi)\geq \gamma_M \sum_{i=1}^q|X\wedge e_i|^2+|hX|^2>{\rm cst}'|X|^2>0,$$ 
for all $X\in \Gamma(Q)$ which means that $H^1(M)=0.$ Using the first result in the Appendix, we find a contradiction. Thus, we deduce that $|b_1|=\cdots=|b_m|=0$ which means that the normal bundle is integrable. In this case, the universal cover of $M$ is isometric to the Riemannian product of $\mathbb{R}\times \Sigma$ where $\Sigma$ is a simply connected compact manifold with positive curvature. This ends the proof. 
\hfill$\square$ 

\section{Appendix}\label{appendix} 
The following results are partially contained in \cite[Rem. 2.14]{J}, \cite[Prop. 1.8]{B} and \cite{E} but we include them here for completeness. Let us denote by $b_s(M)={\rm dim}\, H^s(M)$ (resp. $b_s(\mathcal{F})={\rm dim}\, H_b^s(\mathcal{F})$) the betti numbers (resp. basic betti numbers). 

\begin{prop} Let $(M,g,\xi)$ be a compact Riemannian manifold endowed with a Riemannian flow of codimension $q$ with basic mean curvature $\kappa$. Assume that the first cohomology group $H^1(M)=\{0\}.$ Then we have that $b_2(\mathcal{F})=1+b_2(M).$
\end{prop}

{\noindent \it Proof.} We use the long exact sequence of cohomologies stated in \cite[Thm. 3.2]{Pri} 
$$0\to H_b^1(\mathcal{F}) \to H^1(M)\bui{\to}{j} H^{q}_b(\mathcal{F})\bui{\to}{i_1} H^2_b(\mathcal{F})\bui{\to}{i_2} H^2(M)\to H^{q-1}_b(\mathcal{F}),$$ 
where $i_1=\wedge [\Omega]$ and $i_2$ is the inclusion map. Since $H^1(M)=0,$ we have that $H^{q}_b(\mathcal{F})\simeq \R$ and $H^{q-1}_b(\mathcal{F})\simeq H^{1}_b(\mathcal{F})=\{0\}$ (see \cite{T}). From the fact that the map $i_1$ is injective, $i_2$ is surjective and ${\rm Im}\, i_1={\rm Ker}\, i_2,$ we find that ${\rm Ker}\, i_2\simeq \mathbb{R}$ and ${\rm Im}\, i_2=H^2(M).$ Therefore, we deduce the statement of the proposition.
\hfill$\square$  

\begin{prop} Let $(M,g,\xi)$ be a compact Riemannian manifold endowed with a minimal Riemannian flow of codimension $q$. Assume that ${\rm Ric}^M(\xi)=\lambda\xi$ with $\lambda>0.$ Then the Euler class $[d\xi]$ is a non-zero cohomology class in $H_b^2(\mathcal{F}).$ Moreover, we have that $b_1(\mathcal{F})=b_1(M)$ and $1\leq b_2(\mathcal{F})\leq 1+b_2(M).$
\end{prop}

{\noindent \it Proof.} Take an orthonormal frame $\{e_i\}_{i=1,\cdots,q}$ in $\Gamma(Q)$ and consider $Y=Z=e_i$ in the formula $g(R^M(X,Y)\xi,Z)=g(-(\nabla_Xh)Y+(\nabla_Y h)X,Z).$ After tracing over $i,$ we get that ${\rm Ric}^M(\xi,X)=(\delta h)(X)$ for all $X\in \Gamma(Q).$ The assumption ${\rm Ric}^M(\xi)=\lambda\xi$ gives that the basic $2$-form $\Omega:=-\frac{1}{2}d\xi=g(h\cdot,\cdot)$ is co-closed. As $\Omega$ is also a closed form, it then becomes a basic-harmonic form. But the choice of $\lambda=|h|^2$ to be strictly positive implies that the form $\Omega$ does not vanish. This shows the first part. To prove the second part, we use again the Gysin sequence as in the previous proposition and the fact that $H^{q}_b(\mathcal{F})\simeq \R$ (recall the flow is minimal) to get that $i_1$ is injective and thus $b_2(\mathcal{F})=1+{\rm dim}\,{\rm Im}\,i_2$. Also, we get that $j=0$ and therefore $H_b^1(\mathcal{F}) \simeq H^1(M).$ This finishes the proof. 
\hfill$\square$


\begin{thebibliography}{99} 
\bibitem{B} C. Boyer, K. Galicki and M. Nakayame, {\it On positive Sasakian Geometry}, Geom. Dedicata {\bf 101} (2003), 93--102. 
\bibitem{BW} C. B\"ohm and B. Wilking, {\it Manifolds with positive curvature operators are space forms}, Annals of Math. {\bf 167} (2008), 1079--1097. 
\bibitem{CS} Q. Cui and L. Sun, {\it A sharp lower bounds of eigenvalues for differential forms and homology sphere theorems}, arxiv:1704.00668v1.
\bibitem{D} D. Dom\'inguez, {\it A tenseness theorem for Riemannian foliations}, C. R. Acad. Sci. S\'er. I {\bf 320} (1995), 1331--1335. 
\bibitem{E} A. El Kacimi, {\it Op\'erateurs transversalement elliptiques sur un feuilletage riemannien et applications}, Compositio Mathematica {\bf 79} (1990), 57--106. 
\bibitem{GM} S. Gallot and D. Meyer, {\it Op\'erateur de courbure et laplacien des formes diff\'erentielles d'une vari\'et\'e riemanienne}, J. Math. Pures. Appl. {\bf 54} (1975), 259--284. 
\bibitem{GK} J. F. Glazebrook and F.W. Kamber, {\it Transversal Dirac families in Riemannian foliations}, Comm. Math. Phys.  {\bf 140} (1991), 217--240.
\bibitem{GG} D. Gromoll and K. Grove, {\it One dimensional metric foliations in constant curvature space}, Diff. Geom. Comp. Anal. H.E. Rauch memorial volume, Springer, Berlin (1985), 165--167.
\bibitem{Hab} G. Habib, {\it Energy-Momentum tensor on foliations}, J. Geom. Phys. {\bf 57} (2007), 2234-2248.
\bibitem{HR} G. Habib and K. Richardson, {\it Modified differentials and basic cohomology for Riemannian foliations}, J. Geom. Anal. {\bf 23} (2013), 1314--1342.  
\bibitem{L} A. Lichnerowicz, {\it G\'eom\'etrie des groupes de transformations}, Travaux et Recherches Math\'ematiques, III, Dunod, Paris, 1958. 
\bibitem{RP} E. Park and K. Richardson, {\it The basic Laplacian of a Riemannian foliation}, Amer. J. Math. {\bf 118} (1996), 1249--1275.
\bibitem{J} P. Jammes, {\it Effondrement, spectre et propri\'et\'es diophantiennes des flots riemanniens}, Ann. Inst. Fourier. {\bf 60} (2010), 257--290.
\bibitem{JR} S. D. Jung and K. Richardson, {\it Transversal conformal Killing forms and a Gallot-Meyer theorem for foliations}, Math. Z. {\bf 270} (2012), 337--350.
\bibitem{M} A. Mason, {\it An application of stochastic flows of Riemannian foliations}, Houst. J. Math. {\bf 26} (2000), 481--515. 
\bibitem{ORT} M. Min-Oo, E. Ruh and Ph. Tondeur, {\it Vanishing theorems for the basic cohomology of Riemannian foliations}, J. Reine Angew. Math.   {\bf 415} (1991), 167--174. 
\bibitem{Ob} M. Obata, {\it Certain conditions for a Riemannian manifold to be isometric with a sphere}, J. Math. Soc. Japan  {\bf 14} (1962), 333--340.
\bibitem{O} B. O'Neill, {\it The fundamental equations of a submersion}, Mich. Math. J.  {\bf 13} (1966), 459--469. 
\bibitem{P} P. Petersen, {\it Riemannian Geometry}, Graduate Texts in Mathematics {\bf 171}, Springer, New York (1998). 
\bibitem{Pri} J. I. Royo Prieto, {\it The Gysin sequence for Riemannian flows}, Contem. Math. {\bf 288} (2001), 415--419. 
\bibitem{S} A. Savo, {\it The Bochner formula for isometric immersions}, Pacific J. Math.  {\bf 272} No. 2 (2014), 395--422.
\bibitem{Se} U. Semmelmann, {\it On conformal Killing tensor in a Riemannian space}, Math. Z. {\bf 245} (2003), 503--527.
\bibitem{Ta} S. Tachibana, {\it On Killing tensors in Riemannian manifolds with positive curvature operator}, Tohoku. Math. J. {\bf 28} (1976), 177--184.
\bibitem{T} Ph. Tondeur, {\it Geometry of Foliations}, Birkh\"auser, Boston, 1997. 
\end{thebibliography}
\end{document}